\documentclass[11pt]{article}
\usepackage{amssymb,amsmath,graphicx,amsthm,mathrsfs}
\newtheorem{Theorem}{Theorem}[section]
\newtheorem{Lemma}[Theorem]{Lemma}

\newtheorem{Remark}[Theorem]{Remark}
\numberwithin{equation}{section}

 \def\p{\partial} 
\def \Vh0{\stackrel{\circ}{V}_h}

\newcommand{\q}{\quad}   \newcommand{\qq}{\qquad} 

  \def\a{\alpha} 
   
\def\ms{\medskip}  \def\ss{\smallskip}

\newcommand{\lc}
{\mathrel{\raise2pt\hbox{${\mathop<\limits_{\raise1pt\hbox
{\mbox{$\sim$}}}}$}}}

\newcommand{\gc}
{\mathrel{\raise2pt\hbox{${\mathop>\limits_{\raise1pt\hbox{\mbox{$\sim$}}}}$}}}

\newcommand{\ec}
{\mathrel{\raise2pt\hbox{${\mathop=\limits_{\raise1pt\hbox{\mbox{$\sim$}}}}$}}}

\def\cU{{\cal U}}
\def\cV{{\cal V}}

\def\cP{{\cal P}}

\def\bb{\begin{equation}} \def\ee{\end{equation}}

\def\beqn{\begin{eqnarray}}  \def\eqn{\end{eqnarray}}

\def\beqnx{\begin{eqnarray*}} \def\eqnx{\end{eqnarray*}}

\def\bn{\begin{enumerate}} \def\en{\end{enumerate}}

\def\bd{\begin{description}} \def\ed{\end{description}}

\def \om {\Omega}

\def \p  {\partial}
\def \l  {\lambda}
\def \d  {\displaystyle}

\def \e  {\varepsilon}

\def\cS{{\cal S}}

\def\lan{\mathop{\langle}}
\def\ran{\mathop{\rangle}}

\def\bp{\begin{proposition}}
\def\ep{\end{proposition}}
\def\ba{\begin{array}}
\def\ea{\end{array}}

\def\sqr#1#2{{\vcenter{\vbox{\hrule height.#2pt
              \hbox{\vrule width.#2pt height#1pt \kern#1pt \vrule width.#2pt}
              \hrule height.#2pt}}}}
\def\signed #1{{\unskip\nobreak\hfil\penalty50
              \hskip2em\hbox{}\nobreak\hfil#1
              \parfillskip=0pt \finalhyphendemerits=0 \par}}
\def\endpf{\signed {$\sqr69$}}
\def\see{{\it see} }

\def\ds{\displaystyle}
\def\nm{\noalign{\ms}}
\def\ns{\noalign{\ss}}
\def\d{\delta}
\def\l{\lambda}

\begin{document}

\title{\bf Advantages for controls imposed in a proper subset\footnote{This study is partially
supported  by the National Natural Science Foundation of China under grants
10871154, 10831007, 10801041, 11161130003  and 11171264; and by the National Basis Research
 Program of China (973 Program) under grant 2011CB808002.}}

\author{Gengsheng Wang \footnote{ School of Mathematics and Statistics, Wuhan University, Wuhan,
430072, China. E-mail: wanggs62@yeah.net  },\qq Yashan Xu\footnote{  School of Mathematical
Sciences, Fudan University, Shanghai 200433, China. E-mail: yashanxu@fudan.edu.cn}\date{}}

\maketitle

\begin{abstract}  In this paper, we study  time optimal control problems for heat equations on $\Omega\times \mathbb{R}^+$. Two properties under consideration are
the existence and the bang-bang properties of time optimal controls. It is proved that those two properties hold when controls are imposed on some proper subsets of $\Omega$; while they do not stand
when controls are active on the whole $\Omega$. Besides, a new property for eigenfunctions of $-\Delta$ with  Dirichlet boundary condition is revealed.

\end{abstract}

\ms


\ms


\rm

\section{Introduction.}\label{1}
$\;\;\;\;\;$Let $\Omega$ be a bounded domain in $\mathbb{R}^n$. Let $\omega$ be a non empty and open subset of
  $\Omega$. Write $\chi_{\omega}$ for its  characteristic function. Consider the following controlled heat equation:
\begin{equation}\label{state} \left\{\begin{array}{lll}
\ns \p_ty(x,t)-\triangle y(x,t)=\chi_\omega(x)u(x,t)&\mbox{in}&\Omega\times\mathbb{R}^+,\\
\ns y(x,t)=0&\mbox{on}&\partial \Omega\times \mathbb{R}^+,\\
\ns y(x,0)=y_0(x) &\mbox{in}&\om,
\end{array}\right.\end{equation}
where $u$ is a  control function  taken from a  control constraint set and $y_0$ is an initial state taken from  $L^2(\Omega)$.
The solution of (\ref{state}) corresponding to $u$ and $y_0$ will be treated as a function from $\mathbb{R}^+$ to $L^2(\Omega)$ and denoted by $y(\cdot; u,y_0)$.

The purpose of this study is to reveal  the following fact: {\it Some properties hold for some time optimal control problems of (\ref{state}) when $\omega$ is a proper subset of $\Omega$,
 but do not stand when $\omega=\Omega$. Consequently,  the local  control may  be more effective  than the global control for heat equations in some cases}.

We begin with  introducing  time optimal control problems. Let  $\bigr\{\xi_i \bigl\}^\infty_{i=1}$ be a complete set
of eigenfunctions for $-\Delta$ with  Dirichlet boundary condition such that it serves as a normalized  orthonormal basis of $L^2(\Omega)$. Write $\bigr\{\lambda_i \bigl\}^\infty_{i=1}$, with $ 0<\l_1 <\l_2 \le\cdots<+\infty$,
for the corresponding set of  eigenvalues. Then, we  take  the following target set:
$$\cS_m=span \bigr\{\,\xi_{m+1},\xi_{m+2},\cdots\bigl\},\;\;\mbox{where}\;\; m\geq 2\;\;\mbox{is arbitrarily  fixed}.$$
\noindent Next, we define, for each natural number $k$ and each finite sequence  of positive  numbers $\{\bar a_i\}_{i=1}^k$,
 the following  control constraint  set:
$$\cU_{\{\bar a_i\}_{i=1}^k}=\left\{\sum\limits_{i=1}^k  \alpha_i(\cdot)\xi_i \; \biggm| \;\mbox{each}\;\; \alpha_i(\cdot) \;\;\mbox{is measurable from}\;\;\mathbb{R}^+\;\;\mbox{to}\;\;
\left[-\bar a_i, \bar a_i\right] \right\}.$$
\noindent Consider the following time optimal control problem:
 $$
 {(\cP)}\hspace{1cm}\ds\inf\limits\left\{t\geq 0 \; \big| \; y(t; u,y_0)\in \cS_m\right\},\;\;\mbox{where the infimum is taken over all}\;\; u\in\cU_{\{\bar a_i\}_{i=1}^k}.
 $$
 \noindent   Two properties of Problem $(\mathcal{P})$
  under consideration  are as follows: $(i)$ The existence of time optimal controls; $(ii)$ The bang-bang property: {\it  any optimal control $u^*=\sum_{i=1}^k \alpha^*_i\xi_i$  satisfies that
  for each $i$, $|\alpha^*_i(t)|=\bar a_i$ for almost every $t\in (0,t^*)$, where $t^*$ is the optimal time.} In the case that $\Omega$, $\omega$,  $k$ and $y_0 \notin \mathcal{S}$  are fixed, we say Problem
  $(\mathcal{P})$ has optimal controls if for any finite sequence of positive numbers  $\{\bar a_i\}_{i=1}^k$, it has optimal controls. When $\omega=\Omega$, $y_0\notin\mathcal{S}_m$ and $k$ are given,  Problem
  $(\mathcal{P})$ has optimal controls if and only if there is a finite sequence of positive numbers  $\{\bar b_i\}_{i=1}^k$ such that the problem $(\mathcal{P})$, with $\{\bar b_i\}_{i=1}^k$, has optimal controls (see Remark~\ref{WGremark2.3}).

The main results of this paper are broadly stated as follows: $(a)$ Suppose that  $\omega=\Omega$ and $y_0\notin \mathcal{S}_m$. Then, $k$ and $y_0$ are such that  Problem $(\cP)$ has no optimal control if and only if  $k<m$ and $y_0$ satisfies
\begin{eqnarray}\label{w2.1}\left(\lan y_0,\xi_{k+1}\ran,\lan y_0,\xi_{k+2}\ran,\cdots,\lan y_0,\xi_{m}\ran\right)^T \neq0;
\end{eqnarray} $(b)$ Suppose that  $\omega=\Omega$ and $y_0\notin \mathcal{S}_m$. Assume that either $k\geq m$ or $k<m$  and   $y_0$ does not satisfy (\ref{w2.1}). Then, in general, Problem $(\mathcal{P})$ does not hold the bang-bang property; $(c)$ Suppose that $\Omega$ and  $\omega$  satisfy accordingly the following conditions:
\begin{itemize}
\item {\bf (D1)} ~The  eigenvalues $\l_1\cdots\l_m$ are simple, i.e., $\l_1<\l_2<\cdots<\l_m$,
 \end{itemize}
 and
\begin{itemize}
\item {\bf (D2)}
 ~$ \lan \chi_{\omega}\xi_i\,,~\xi_j\ran\neq 0\q{\rm for~ all ~} i\in \{1,2,\cdots,m\}\;\;\mbox{and}\;\;j\in\{1,2,\cdots,k\}.$
 \end{itemize}
\noindent Then, for each $k\geq 1$, each $y_0\notin \mathcal{S}_m$ and each finite sequence of positive numbers $\{\bar a_i\}_{i=1}^k$,  Problem $(\mathcal{P})$ has optimal controls and holds the bang-bang property.\\

It is worth mentioning that for any fixed bounded domain $\Omega$, there are a lot of open subsets $\omega$ in $\Omega$
such that $ \lan \chi_{\omega}\xi_i\,,~\xi_j\ran\neq 0$ for all $i,j=1,2,\cdots$ (see Theorem~\ref{theorem4.1} for a new property of the eigenfunctions $\{\xi_i\}_{i=1}^\infty$); while there are a lot of bounded domains $\Omega$ such that the property $\bf(D1)$ holds (see Remark~\ref{remark4.2})).

\section {Studies of Problem $(\mathcal{P})$ where $\Omega=\omega$}
$\;\;\;\;\;$ The following result is another version of Theorem 2.5 in  \cite{Evans}. It will be used later.
\begin{Lemma}\label{lemma3.1}  Let   $\hat A\in \mathbb{R}^{d \times d}$ and
 $\hat B\in \mathbb{R}^{d\times l}$, where $d$ and $l$ are  natural numbers.
  Suppose that
\begin{equation}\label{3-2}
 rank\left(\hat B,\hat  A\hat B,\hat  A^2\hat B, \cdots,\hat  A^{d-1}\hat B\right)=d,
 \end{equation}
 and the spectrum of  $\hat A$ belongs to the left half plane of $\mathbb{C}$. Then, for each finite sequence of positive numbers $\{b_i\}_{i=1}^l$ and
  each   $w_0$ in $\mathbb{R}^d$, there are a  $\hat t\geq 0$ and  a control $\hat \beta$ in the set:
\begin{equation}\label{3-1}
\bar \cV \triangleq\left\{ \beta=\left( \beta_1,\cdots,  \beta_{l}\right)^T\bigm| \;\mbox{each}\;\; \beta_i \; \mbox{is measurable from}\; \mathbb{R}^+\;\;\mbox{to}\;\;
[-b_i,b_i]\right\},\end{equation}
 such that the solution $\hat w(\cdot;\hat\beta,w_0)$ to the equation:
\begin{equation}\label{3-3}
 \left\{\begin{array}{ll}
\ns\dot{\hat w}(t)=\hat A\hat w(t)+\hat B\hat \beta(t),\q&t\in\mathbb{R}^+,\\
\ns \hat  w(0)=w_0,
\end{array}\right.\end{equation}
reaches zero at $\hat t$.

 \end{Lemma}

\begin{Theorem}\label{theorem2.1} Suppose $\omega=\Omega$ and  $y_0\notin \cS_m$. Then, $k$ and $y_0$ are such that
Problem $(\cP)$ has no optimal control if and only if  $k<m$ and $y_0$ satisfies (\ref{w2.1}).
\end{Theorem}
\noindent{\it Proof.} The proof will be organized in three steps as follows:\\

\noindent{\it Step 1. Suppose that  $k<m$ and $y_0$ satisfies (\ref{w2.1}). Then, for any finite sequence of positive numbers $\{\bar a_i\}_{i=1}^k$,  Problem $(\cP)$ has no optimal control.}\\

Let  $\{\bar a_i\}_{i=1}^k$ be a finite sequence of positive numbers.
Then each $u(\cdot)\in \cU_{\{\bar a_i\}_{i=1}^k}$ can be expressed as  $u(t)=\sum\limits_{i=1}^k  \alpha_i(t)\xi_i$. Write $y(t ;u, y_0)=\sum\limits^\infty_{i=1}y_i(t)\xi_i$.
Clearly, the controlled equation (\ref{state}) is equivalent to the following system:
$$\ds\dot{y}_i(t)+\l_i y_i(t)
=\sum\limits_{j=1}^k \alpha_j(t)\lan\chi_\omega\xi_i,\xi_j\ran,\; \; y_i(0)=\left<y_0,\xi_i\right>,\qq i=1,2,\cdots.
$$
Write
$$
z(t)=\left(\ba{c}y_1(t)\\y_2(t)\\\cdots\\y_m(t)\ea\right),\q
A=\left(\ba{cccc}\l_1\\&\l_2\\&&\ddots\\&&&\l_m\ea\right),\q
\a(t)=\left(\ba{c}\alpha_1(t)\\\alpha_2(t)\\\cdots\\\alpha_k(t)\ea\right),
$$
and
\begin{equation}\label{1.2}
B=\Bigr(\lan \chi_\omega \xi_i\,,~\xi_j\ran \Bigl)_{i,j}\in\mathbb{R}^{m\times k}.
\end{equation}
Let
$$U_{\{\bar a_i\}^k_{i=1}}=[-\bar a_1,\bar a_1]\times[-\bar a_2,\bar a_2]\times\cdots\times[-\bar a_k,\bar a_k].$$
Consider  the following time optimal control problem:
$$(\widetilde{\cP})
\hspace{1.5cm}\ds\inf\limits\left\{t\geq 0\; \big| \;  z(t;\a,z_0)=0\right\},\hspace{3.5cm}$$
where the infimum is taken over all $\alpha$ from the  control constraint set:
\begin{eqnarray*}
\cV_{\{\bar a_i\}_{i=1}^k}\triangleq\left\{\alpha=(\alpha_1,\cdots, \alpha_k)^T\bigm| \;\mbox{each}\;\; \alpha_i \; \mbox{is measurable from}\; \mathbb{R}^+\;\;\mbox{to}\;\;
 [-\bar a_i, \bar a_i]\right\},\end{eqnarray*}
and $z(\cdot; \alpha,z_0)$ is the solution to the following  equation:
\begin{equation}\label{statef}
 \left\{\begin{array}{lll}
\ns\dot{z}(t)+Az(t)=B\a(t),\;\;t\in\mathbb{R}^+,\\
\nm z(0)= \left(\lan y_0,\xi_1\ran, \cdots,\lan y_0, \xi_m\ran\right)^T.
\end{array}\right.\end{equation}
Clearly, Problems $(\mathcal{P})$ and $(\widetilde{\cP})$ are equivalent, i.e., $t^*$ and $u^*=\sum_{i=1}^k \alpha^*_i\xi_i$  are accordingly the optimal time and an optimal control to
Problem $(\mathcal{P})$ if and only if $t^*$ and $(\alpha^*_1,\cdots,\alpha^*_k)^T$ are the optimal time and an optimal control to Problem $(\widetilde{\cP})$ respectively.

Since $\omega=\Omega$ and $k<m$, it follows from (\ref{1.2}) that
$B=\left(\ba{c}I_{k\times k}\\0\ea\right)$ in this case.
Let $z_1(t)=(y_1(t),\cdots, y_k(t))^T$ and $z_2(t)=(y_{k+1}(t),\cdots, y_m(t))^T$. Write
$$A_1=\left(\ba{cccc}\l_1\\&&\ddots\\&&&\l_k\ea\right)\;\;\mbox{and}\;\; \;\; A_2=\left(\ba{cccc}\l_{k+1}\\&&\ddots\\&&&\l_m\ea\right).
$$
Then, Equation (\ref{statef}) can be written as
\begin{equation}\label{0.1}
\ds\frac{d}{dt}\left(\ba{c}z_1\\z_2\ea\right)(t)+\left(\ba{cc}A_1\\&A_2\ea\right)\left(\ba{c}z_1\\z_2\ea\right)(t)
=\left(\ba{c}I_{k\times k}\\0\ea\right)\a(t),\end{equation}
together with the initial condition:
$$ (z_1(0), z_2(0))^T= \left( \left(\lan y_0,\xi_1\ran, \cdots,\lan y_0, \xi_k\ran\right)^T,  \left(\lan y_0,\xi_{k+1}\ran, \cdots,\lan y_0, \xi_m\ran\right)^T\right)^T
$$
 This, along with the condition (\ref{w2.1}), indicates that\
$z_2(t)\neq 0$,  for each $t>0$ and each control $\alpha$ in $\cV_{\{\bar a_i\}_{i=1}^k}$. Consequently, Problem $(\mathcal{P})$ has no time optimal control. \\

\noindent{\it Step 2. Suppose that $k<m$ and $y_0$ does not satisfy (\ref{w2.1}). Then,   Problem   $(\mathcal{P})$ has
 optimal controls. }\\

Let  $\{\bar a_i\}_{i=1}^k$ be a finite sequence of positive numbers.  Since $y_0$ does not satisfy (\ref{w2.1}), it holds that  $z_2(0)=0$. Thus, it follows from (\ref{0.1}) that  $z_2(t)=0$ for all $t\geq 0$. Hence, Problem  $(\widetilde \cP)$ shares the same optimal time and optimal controls with
 the following time optimal control problem:
$$(\widetilde{\cP}_1):\hspace{0.5cm}\ds \inf\limits\left\{t\geq 0\; \big| \;  z_1(t;\a)=0\right\},$$
where the infimum is taken over all $\alpha$ from $\cV_{\{\bar a_i\}_{i=1}^k}$,
and $z_1(\cdot; \alpha)$ is the solution to the   equation:
$$ \dot{z}_1(t)+A_1z_1(t)=I_{k\times k}\a(t),\;\;t\in\mathbb{R}^+,\;  z_1(0)=\left(\lan y_0,\xi_1\ran, \cdots,\lan y_0, \xi_k\ran\right)^T.$$
According to Lemma~\ref{lemma3.1}, Problem $(\widetilde{\cP}_1)$ has admissible controls. Then,  by the standard argument
(see either Theorem 13 and the note after it in Chapter III on Page 130 in
 \cite{Pontryagin} or   Theorem 3.1 on Page 31 in \cite{Evans}), one can easily verify that Problem $(\widetilde{\cP}_1)$ has optimal controls. Consequently, Problem  $(\mathcal{P})$ has
 optimal controls.
\\

\noindent{\it Step 3. Suppose that  $k\ge m$. Then,  Problem  $(\cP)$ admits optimal controls.}\\

Let  $\{\bar a_i\}_{i=1}^k$ be a finite sequence of positive numbers.
 Since $B=\left(I_{m\times m},0_{m\times (k-m)}\right)$ in the case that $k\geq m$,  control variables
$\alpha_{m+1}(\cdot),\cdots, \alpha_{k}(\cdot)$ play no role in Equation (\ref{statef}) when $k>m$.  Hence, in the case that $k\geq m$, the effective controls in Problem  $(\widetilde{P})$
 have the form:
$\hat\a =(\alpha_1(\cdot),\cdots, \alpha_{m}(\cdot))^T$. Therefore, Problem $(\widetilde{\cP})$ shares the same optimal time and optimal controls with the following time optimal control  problem:
$$(\widetilde{\cP}_2):\hspace{0.5cm}\ds
\inf\limits\left\{t\geq 0\; \big| \;  z(t;\hat\a )=0\right\},$$
where the infimum is taken over all $\hat\alpha\triangleq (\alpha_1,\cdots,\alpha_m)^T$ from the control constraint set:
\begin{eqnarray*}
\cV_{\{\bar a_i\}_{i=1}^m}\triangleq\left\{\alpha=(\alpha_1,\cdots, \alpha_m)^T\bigm| \;\mbox{each}\;\; \alpha_i \; \mbox{is measurable from}\; \mathbb{R}^+\;\;\mbox{to}\;\;
 [-\bar a_i, \bar a_i]\right\},\end{eqnarray*}
 and $z(\cdot; \hat\a)$ is the solution of the following equation:
\begin{equation}\label{2.2}
 \left\{\begin{array}{l}
\ns\dot{z}(t)+Az(t)=I_{m\times m}\hat\a (t),\q t\in\mathbb{R}^+,\\
\nm z(0)=z_0\triangleq \left(\lan y_0,\xi_1\ran, \lan y_0, \xi_2\ran,\cdots,\lan y_0, \xi_m\ran\right)^T.
\end{array}\right.\end{equation}
Then, by Lemma~\ref{lemma3.1}, using the same argument as that in Step 2, one can prove that
 Problem $(\widetilde{\cP})$ has optimal controls. \\

 In summary, we complete the proof.\endpf

\bigskip

\begin{Remark}\label{WGremark2.3}
From the proof of Theorem~\ref{theorem2.1}, one can easily verify the following: $(a)$ Suppose that  $\omega=\Omega$. Let $y_0\notin \mathcal{S}_m$ and $k$ be given. Then, Problem $(\mathcal{P})$ has optimal controls if and only if there is a finite sequence of positive numbers $\{\bar b_i\}_{i=1}^k$ such that the problem $(\mathcal{P})$, with $\{\bar b_i\}_{i=1}^k$, has optimal controls.
$(b)$ In the case that $\omega=\Omega$ and $y_0\notin\mathcal{S}_m$, Problem $(\mathcal{P})$ has optimal controls, provided either $k\geq m$ or $k<m$ and $y_0$ does not satisfy
(\ref{w2.1}).

\end{Remark}
\bigskip

\begin{Theorem}\label{theorem2.2}
Let $\omega=\Omega$ and $y_0\notin \cS_m$. Let  $\{\bar a_i\}_{i=1}^k$ be a finite sequence of positive numbers.   For each $i\in \{1,\cdots, k\}$, write
\begin{equation}\label{d3}
 T_i= \ds\frac {1}{\l_i}\ln\left(1+\frac{\l_i}{\bar a_i}\left|\lan y_0,\,\xi_i\ran\right|\right).\qq
\end{equation}
 Then Problem $(\cP)$, with $\{\bar a_i\}_{i=1}^k$,  does not have the
   bang-bang property,  if either of the following conditions stands: $(i)$ $k\ge m$ and the numbers $T_1,\cdots, T_m$ are not the same; $(ii)$ $k<m$, $y_0$ does not satisfy
(\ref{w2.1}) and the numbers $T_1,\cdots, T_k$ are not the same.
\end{Theorem}
\noindent{\it Proof.} Simply write $(\mathcal{P})$ for the problem  $(\cP)$, with $\{\bar a_i\}_{i=1}^k$.
For each $i\in \{1,\cdots k\}$, define
\begin{eqnarray}\label{w2.6}\widetilde{ \alpha}_i(\cdot)=-\chi_{[0,T_i]}(\cdot) sgn\big(\lan y_0,\xi_i\ran\big)\bar a_i
\triangleq\left\{\ba{cl}\ns\chi_{[0,T_i]}(\cdot)\bar a_i,\q  & \mbox{if}\; \lan y_0,\xi_i\ran<0,\\
\ns0,\q  & \mbox{if}\; \lan y_0,\xi_i\ran=0,\\
\ns-\chi_{[0,T_i]}(\cdot)\bar a_i,\q  & \mbox{if}\; \lan y_0,\xi_i\ran>0.\ea
\right.
\end{eqnarray}
We first prove the following property $(\mathcal{H}_1)$: {\it When $k\geq m$, $\widetilde{T}$ and $\widetilde{u}$ are the optimal time and an optimal control to Problem    $(\mathcal{P})$ respectively,
where
$$
\widetilde{T}\triangleq\max\{T_1,T_2,\cdots, T_m\}\;\;\mbox{and}\;\;\widetilde{u}\triangleq\sum\limits_{i=1}^m\widetilde{ \alpha}_i\xi_i.
$$}
By the equivalence of Problems $(\mathcal{P})$ and $(\widetilde{\cP}_2)$ (see Step 3 in the proof of Theorem~\ref{theorem2.1}), we need only to verify that $\widetilde{T}$ and $\widetilde{\alpha}$ are the optimal time and an optimal control to Problem  $(\widetilde{\cP}_2)$  respectively, where $\widetilde{\alpha}\triangleq (\widetilde{ \alpha}_1,\cdots, \widetilde{ \alpha}_m)^T$.

For this purpose, we observe  from  direct computation that for each $i\in \{1,\cdots, m\}$,
$T_i$ and $\widetilde{\alpha}_i(\cdot)$ are the   optimal time and the optimal control to the following time optimal control  problem:
$$ (P_i): \hspace{0.5cm}\ds \inf\limits\left\{t\geq 0\; \big|\;  z_i(t;\alpha_i)=0\right\},$$
 where the infimum is taken over all $\alpha_i(\cdot)$ from the set of all measurable functions from $R^+$ to $[-\bar a_i, \bar a_i]$, and
$z_i(\cdot; \alpha_i)$ solves the following equation:
$$\dot{z}_i(t)+\l_i(t)z_i(t)=\alpha_i(t),\;\; z_i(0)=\lan y_0, \xi_i\ran.$$
Clearly, $\widetilde{\alpha}\in\cV_{\{\bar a_i\}_{i=1}^m}$ and   $(\left(z_1\left(\cdot ;\widetilde{\a}_1\right),\cdots,z_m\left(\cdot ;\widetilde{\a}_m\right)\right)^T$ is the solution $z\left(\cdot; \widetilde{\alpha}\right)$ to Equation (\ref{2.2}) with $\hat \a=\widetilde{\alpha}$.
  Since $z_i\left(T_i;\widetilde{\alpha}_i\right)=0$, it holds that
  \begin{eqnarray}\label{WG2.6}
  z_i\left(\widetilde{T}; \widetilde{\a}_i\right)=0\;\;\mbox{ for all}\;\; i\in\{1,2,\cdots, m\},\;\;\mbox{i.e.,}\;\;  z\left(\widetilde{T}; \widetilde{u}\right)=0.
  \end{eqnarray}
Hence, the optimal time to Problem $(\widetilde{\cP}_2)$ is not bigger than $\widetilde{T}$.
On the other hand, if $\hat\a\triangleq(\hat\a_1\cdots,\hat\a_m)^T\in \cV_{\{\bar a_i\}_{i=1}^m}$ and $\hat T>0$ are such that
$z\left(\hat T;\hat\a\right)=0$,
then it stands that
$$
z_i\left(\hat T;\hat \alpha_i\right)=0\;\;\mbox{ for all}\;\; i\in \{1,\cdots, m\}.$$
 By the optimality of $T_i$ to Problem $(P_i)$, we see that
$\hat T\ge T_i$ for all $i\in  \{1,\cdots, m\}$, from which, it follows that
$\hat T\ge \widetilde{T}$. Therefore, $\widetilde{T}$ is the optimal time to Problem $(\widetilde{\cP}_2)$. Along with (\ref{WG2.6}), this yields that $\widetilde{\alpha}$ is an optimal control
to this problem. Hence, the property  $(\mathcal{H}_1)$ stands.

 Since $y_0\notin \cS_m$, it holds  that $\widetilde{T}>0$. Because $T_1,\cdots, T_m$ are not the same,  there is an $i_0\in\{1,2,\cdots,m\}$ such that $T_{i_0}<\widetilde{T}$. Then, it follows from (\ref{w2.6}) that
$\widetilde{\alpha}_{i_0}(t)=0$ for all $t\in (T_{i_0},\widetilde{T}]$. Thus, the optimal control $\widetilde{u}$ does not satisfy the bang-bang property.

Using the very similar argument to that in the proof of the property  $(\mathcal{H}_1)$, one can easily show the following property  $(\mathcal{H}_2)$:  {\it When $k<m$, $y_0$ does not satisfy
(\ref{w2.1}), $\hat T$ and $\hat u$ are the optimal time and an optimal control to Problem $(\mathcal{P})$, where $$
\hat{T}\triangleq\max\{T_1,T_2,\cdots, T_k\}\;\;\mbox{and}\;\;\hat{u}\triangleq\sum\limits_{i=1}^k\widetilde{ \alpha}_i\xi_i.
$$}
Then, by the property $(\mathcal{H}_2)$, (\ref{w2.6}) and the assumptions that $y_0\notin \cS_m$ and the numbers $T_1,\cdots, T_k$ are not the same, one can easily show that the optimal control $\hat u$ does not satisfy
the bang-bang property.
This completes the proof.
\endpf

\section {Studies of Problem $(\mathcal{P})$ where $\omega$ is a  proper subset of  $\Omega$}

\bigskip

\begin{Theorem}\label{theorem3.1} Let $\Omega$ satisfy  the condition $(D1)$. Suppose that $\omega$ holds the condition $(D2)$.
Then, for each $k\geq 1$, each  $y_0\notin \mathcal{S}_m$ and each finite sequence of positive numbers  $\{\bar a_i\}_{i=1}^k$,  Problem $(\mathcal{P})$  has optimal controls.
\end{Theorem}
\noindent{\it Proof.}  By the same way as that in Step 1 of  the proof of Theorem~\ref{theorem2.1}, we define the matrices $A$ and $B$, and the problem $(\widetilde{P})$.
Write   $B_{ij}$ for the element in i-th row and j-th column of $B$, namely, $B_{ij}=\lan \chi_{\omega}\xi_i\,,~\xi_j\ran$. Let $B_1=(B_{11},\cdots, B_{m1})^T$.
We first claim that
\begin{equation}\label{3-4}
rank(B_1, AB_1, A^2B_1, \cdots, A^{m-1}B_1)=m.
\end{equation}
In fact, since
$$A^jB_1=\left(\ba{ccc}\l_1\\&\ddots\\&&\l_m\ea\right)^j\left(\ba{c}B_{11}\\\cdots\\B_{m1}\ea\right)
= \left(\ba{c}\l^j_1B_{11}\\ \cdots\\ \l_m^jB_{m1}\ea\right),$$
it holds that
$$\Bigr|(B_1,AB_1,\cdots,A^{m-1}B_1)\Bigl|
=\left|\ba{cccc}B_{11}&\l_1B_{11}&\cdots&\l_1^{m-1}B_{11}\\
B_{21}&\l_2B_{21}&\cdots&\l_2^{m-1}B_{21}\\
\cdots&\cdots&\cdots&\cdots\\
B_{m1}&\l_mB_{m1}&\cdots&\l_m^{m-1}B_{m1}\ea\right|,$$
which is a   determinant of Vandermonde' type and  equals to
$\prod_{i=1}^{m}B_{i1}\prod_{k>l}(\l_k-\l_l)$. Because of conditions $(D1)$ and $(D2)$, this determinant is not zero, which implies (\ref{3-4}).

Now, according to   Lemma \ref{lemma3.1}, Problem $(\widetilde{\mathcal{P}})$ has  admissible controls.
Then,  by the standard argument (\see either Theorem 13 and the note after it in Chapter III on Page 130 in
 \cite{Pontryagin} or   Theorem 3.1 on Page 31 in \cite{Evans}), one can easily show that Problem $(\widetilde{\cP})$ admits optimal controls.
This, along with  the equivalence of Problems $(\mathcal{P})$ and $(\widetilde{\cP})$, completes the proof.
\endpf

\bigskip
\begin{Remark} From the proof of the above theorem,it follows that Theorem~\ref{theorem3.1} still stands when the condition $(D2)$ is replaced by the following condition:
\begin{itemize}
\item {\bf $(\widetilde{D}2)$}
 ~$ \lan \chi_{\omega}\xi_i\,,~\xi_1\ran\neq 0\q{\rm for~ all } \;\;i\in \{1,2,\cdots,m\}.$
 \end{itemize}
 \end{Remark}

\bigskip
Before studying the bang-bang property for Problem $(\mathcal{P})$ where $\omega$ is a proper open subset of $\Omega$,
we  recall  the general position condition which plays an important role in the studies of the bang-bang  property for linear controlled ordinary differential equations.
Let  $\hat A$ and  $\hat B$ be $m\times m$ and $m\times k$ matrices respectively.  Let $\hat V$ be a  closed polyhedron   in $\mathbb{R}^k$.
 {\it  We say that  $\hat V$  satisfies the general position condition
 with respect to $(\hat A,\hat B)$,  if  for each  nonzero vector $v$, which is parallel to one of  the edges of  $\hat V$,
the vectors
 $$\hat Bv,~\hat A\hat Bv,~\cdots~\hat A^{m-1}\hat Bv$$
are linearly independent.} Consider the following time optimal control
problem:
$$(\hat P): \hspace{2cm}\ds\inf\limits\left\{t: z(t; v,z_0)=0\right\}, \hspace{3.5cm}
$$ where the infimum is taken over all measurable functions $v$
from $\mathbb{R}^+$ to the polyhedron $\hat V$, and $z(\cdot;v,z_0)$ is the solution to the following equation:
$$\dot {z}(t)+\hat Az(t)=\hat Bv(t),\;\;t>0;\;\; z(0)=z_0,$$ with $z_0$ a non-zero vector in $\mathbb{R}^m$.

\begin{Lemma} \label{l2.3}(\see \cite{Pontryagin}, \cite{Fattorini})
Suppose that the closed polyhedron $\hat V$ satisfies the  general position condition  with respect to $(\hat A,\hat B)$. Then any optimal control $\bar u(t)$ to
Problem $(\hat P)$, if exists, takes values on the vertices of $\hat V$ and has a finite number of switchings.
\end{Lemma}

\begin{Theorem}\label{theorem3.3}  Let $\Omega$ satisfy  the condition $(D1)$. Suppose that $\omega$ satisfies the condition $(D2)$.
Then,  for each    $k\geq 1$, each  $y_0\notin \mathcal{S}_m$ and each finite sequence of positive numbers  $\{\bar a_i\}_{i=1}^k$ , Problem $(\mathcal{P})$  holds the   bang-bang property.
\end{Theorem}

\noindent {\it Proof.} By the same way as that in Step 1 of  the proof of Theorem~\ref{theorem2.1}, we define the matrices $A$ and $B$, and the problem $(\widetilde{\mathcal{P}})$.
  According to Lemma \ref{l2.3}, Theorem \ref{theorem3.1} and the equivalence of Problems $(\mathcal{P})$ and $(\widetilde{\mathcal{P}})$,  it suffices to  prove the general condition of  $U_{\{\bar a_i\}^k_{i=1}}$ with respect to $(A, B)$.
 Clearly, the later is equivalent to the statement that for each $j\in \{1,\cdots,k\}$, the vectors $Be_j, ABe_j,\cdots,A^{m-1}Be_j$ are linearly independent, where $\{e_1,\cdots,e_k\}$ is the standard basis of $\mathbb{R}^k$.

Let $F_j=(Be_j, ABe_j,\cdots,A^{m-1}Be_j)$. It is clear that
$$\ba{rl}
\ns|F_j|=&|(Be_j,ABe_j,\cdots,A^{m-1}Be_j)|\\
\nm=&\left|\ba{cccc}B_{1j}&\l_1B_{1j}&\cdots&\l_1^{m-1}B_{1j}\\
B_{2j}&\l_2B_{2j}&\cdots&\l_2^{m-1}B_{2j}\\
\cdots&\cdots&\cdots&\cdots\\
B_{mj}&\l_mB_{mj}&\cdots&\l_m^{m-1}B_{mj}\ea\right|\\
\nm=&\prod_{i=1}^{m}B_{ij}\prod_{k>l}(\l_k-\l_l).\ea$$
This, together with the conditions (D1)  and (D2), yields that $|F_j|\neq 0$ for each $j\in \{1,\cdots,k\}$.
Hence, $U_{\{\bar a_i\}^k_{i=1}}$ satisfies the general position condition  with respect to $(A, B)$. This completes the proof.
\endpf

 \section {Further studies  on the conditions $(D1)$ and $(D2)$}

  $\;\;\;\;\;$In this section,  we first give a remark and a theorem, which reveal accordingly some  properties for eigenvalues and eigenfunctions  of $-\Delta$
  with Dirichlet boundary condition. From the remark, it follows that there are a lot of $\Omega$ satisfying the property $\bf(D1)$. From the theorem,
  it follows that for any  bounded domain $\Omega$ in $\mathbb{R}^n$, there are a lot of $\omega\subset\Omega$
  where the  property $\bf(D2)$ holds.  We end this section with another remark which
  provides an open
  problem.

\medskip

\begin{Remark}\label{remark4.2} It is presented in \cite{Micheletti} (see also \cite{Uhlenbeck}, \cite{Ortega}) that there are a lot of $\Omega$ of class $C^3$ satisfies the condition $(D1)$ in the following sense:
 Let $\Omega$ be a bounded open set of class $C^3$ in $\mathbb{R}^n$. For each $\varepsilon\in (0,1)$,
an $\varepsilon-$neighborhood of  $\Omega$ is defined to be the image $(I+\psi)(\Omega)$, where
   $I$ is the identity map over $\mathbb{R}^n$ and
 $\psi\in C^3(\mathbb{R}^n; \mathbb{R}^n)$, with  the $C^3-$norm  less than $\varepsilon$. For each bounded open set $\widetilde{\Omega}$ of class $C^3$ in $\mathbb{R}^n$,
 Write $\Delta_{\widetilde{\Omega}}$ for the self-adjoint operator in $L^2(\widetilde{\Omega})$ generated by the Laplacian on $\widetilde{\Omega}$  with the homogeneous
 Dirichlet boundary condition. Then, for each $\varepsilon\in (0,1)$, there is an  $\varepsilon-$neighborhood  of $\Omega^\e$ such that $-\Delta_{\Omega^\e}$ has only simple eigenvalues.
\end{Remark}

\medskip

      Before presenting the theorem,
  we introduce the following notations: for each $x\in \mathbb{R}^n$ and each $\rho>0$, $B_\rho(x)$ stands for the open ball in $\mathbb{R}^n$, centered at $x$ and of radius $\rho$;
  $\overline{ B_\rho(x)}$ denotes the closure of the ball $B_\rho(x)$;
   for each $\rho>0$, %
$$\om^\rho\triangleq\Bigr\{x\in \om\setminus\overline{\omega}\bigm| dist \left(\p{ B_\rho(x)},\p \om \right)>0,~ dist \left(\p{  B_\rho(x)},\p \omega \right)>0\Bigl\},$$
where  $dist \left(E_1,E_2 \right)\triangleq\inf\limits_{x_1\in E_1,x_2\in E_2}\|x_1-x_2\|_{\mathbb{R}^n}$ for any subsets $E_1$ and $E_2$ in $\mathbb{R}^n$.
\medskip

\begin{Theorem}\label{theorem4.1}   Let $\Omega$ be a bounded domain in $\mathbb{R}^n$ and $\omega$ be an open subset of $\om$ such that $\om\setminus\overline{\omega}\neq \varnothing$. Then,
for any $\e>0$, there exists an $\varepsilon_0\in (0,\varepsilon)$ such that $\om^{\varepsilon_0}\neq\varnothing$ and for almost every $\widetilde{x}\in \om^{\varepsilon_0}$,
\begin{eqnarray}\label{W4.1}\ds\Bigl\langle \chi_{\omega\cup B_{\varepsilon_0}(\widetilde{x})}\xi_i\,,~ \xi_j\Bigr\rangle\neq0\q{\rm for~ all ~} i,j\in \mathbb{N}.\end{eqnarray}
\end{Theorem}

\noindent {\it Proof.} We  recall that each eigenfunction $\xi_i$ belongs to $ C^\infty(\om)$ (\see Page 335 in \cite{Evans1}). Let
\begin{eqnarray*}
\varphi(x,\tau)=\xi_i(x)\ds e^{\sqrt{\l_i}\tau},\qq (x,\tau)\in \om\times\mathbb{R}.
\end{eqnarray*}
It is obvious that
$$\triangle_x\varphi(x,\tau)+\partial^2_\tau \varphi(x,\tau)=0,\qq (x,\tau)\in \om\times\mathbb{R}.$$
 By the property of harmonic functions (\see Page 6 in \cite{Lin}),
the function  $\varphi(\cdot,\cdot)$ is real analytic over $\om\times\mathbb{R}$.
Thus, each eigenfunction  $\xi_i$ is  real analytic over $\om$.
Write
$$D=\left\{(x,\rho)\in\om\times\mathbb{R}^+\bigm|\overline{ B_\rho(x)} \subset\om\right\}.$$
 Then, for each pair $(i, j)$, we define a function $F_{i,j}(\cdot.\cdot)$ from $D$ to $\mathbb{R}$ by setting:
\begin{equation}
F_{i,j}(x,\rho)=\ds\int_{B_1(0)}\xi_i(x+\rho \eta)\xi_j(x+\rho \eta)d\eta,\;\;(x,\rho)\in D.
\end{equation}
Clearly, it is well defined.  The rest of the  proof will be carried out by the following three steps:

\medskip

\noindent{\it  Step 1.  Suppose that $f$ is a real analytic function  over $\om$. Define the function $F: D\mapsto\mathbb{R}$  by
\begin{equation}
F (x,\rho)=\ds\int_{B_1(0)}f(x+\rho \eta)d\eta,\;\; (x,\rho)\in D.
\end{equation}
 Then $F$ is  real  analytic  over $D$.}

\medskip

  We need only to explain that $F$ is real analytic in a small neighborhood of $(x_0,\rho_0)$ for any point $(x_0,\rho_0)\in D$. First, there is a neighborhood $U$ of $(x_0,\rho_0)$ in $\mathbb{R}^n\times \mathbb{ R}^+$ such that $\overline{  B_\rho (x)}\subset \Omega$ for any $(x,\rho)\in U$.
 Hence, the   function $f(x+\rho \eta)$ is real analytic in $(x,\rho,\eta)$ over $U\times \-{B_1(0)}$.  Extend $f$ to a complex-valued function in $(z,w,\eta)$  over a small neighborhood $U_c\times{B_1(0)}$ of $U\times {B_1(0)}$ in
  ${\mathbb{ C}}^n\times {\mathbb{ C}}\times \-{B_1(0)}$ by making use of  the power series expansion. We then get
$f_c(z+w\eta)$,
%
which is real analytic over $(z,w,\eta)\in U_c\times \-{B_1(0)}$ and holomorphic in $(z,w)\in U_c$ for each fixed  $\eta\in\overline{{B}_1(0)}$.
Clearly, it holds that
$$f_c(x+\rho \eta)=f(x+\rho \eta)\qq{\rm for~ all~ }(x,\rho,\eta)\in U\times B_1(0).$$
Now we define a function $F_c: U_c\mapsto \mathbb{C}$ by setting:
$$F_c (z,w)=\ds\int_{B_1(0)}f(z+w \eta)d\eta,\;\;(z,w)\in U_c,$$
and  define the operator $\bar\p$ in the standard way:
$$\ds\bar\p u(z,w)=\sum^n_{j=1}\frac{\p u}{\p \bar z_j}d\bar z_j+ \frac{\p u}{\p \bar w}d\bar w,$$
where
$$\ds \frac{\p }{\p \bar z_j}=\frac{1}{2}\left( \frac{\p}{\p( Re(z_j))}+\sqrt{-1}\frac{\p}{\p (Im(z_j))}\right),
\qq \frac{\p }{\p \bar w}=\frac{1}{2}\left( \frac{\p}{\p( Re(w))}+\sqrt{-1}\frac{\p}{\p (Im(w))}\right)$$
are the standard Cauchy-Riemann operators (see \cite{Krantz}).
It follows from the holomorphic property  of $f_c$  in $(z,w)$ that
$$\bar\p  F_c(z,w)=\ds\int_{B_1(0)}\bar\p f_c(z+w\eta)d\eta=\int_{B_1(0)}0~d\eta=0.$$
Hence,  $F_c$ is holomorphic in $U_c$. In particular, the function
$F_c(\cdot,\cdot)\Bigl|_{U_c\bigcap(\mathbb{R}^{n}\times\mathbb{R})}$
is real analytic, i.e.  $F(\cdot,\cdot)$ is analytic in $U$.  Thus,  $F(\cdot,\cdot)$ is real analytic over $D$. Consequently,  for each pair $(i,j)$,
the function  $F_{i,j}(\cdot,\cdot)$ is real analytic over $D$.

 \medskip

\medskip

\noindent{\it  Step 2. For each pair $(i,j)$,   $F_{i,j}(\cdot,\cdot)$ is not identically a constant over $D$.} \\

By the unique continuation property of the eigenfunctions  (\see \cite{FHL}), we see that for each $i\in \{1,2,\cdots\}$,
$$ \xi_i(x)\neq 0 \q \mbox{for almost every} \q x\in \Omega.$$
Thus, it holds that for each pair $(i,j)$,  $$ (\xi_i\xi_j)(x)\neq 0 \q \mbox{for almost every} \q x\in \Omega.$$
Since the function
$(\xi_i \xi_j)(\cdot)$ is continuous in $\om$ and $\om\setminus\overline{\omega}\neq \varnothing$,  there is an $\hat x\in \om\setminus\overline{\omega}$
such that $( \xi_i\xi_j)(\hat x)\neq 0 $.
Hence, when $\d>0$ is small enough,  the    function $\bigr(\xi_i \xi_j\bigl)(\cdot)$ is either positive or negative over $B_{\d}(\hat x)$, and $\overline{ B_{\d}(\hat x)}\subset \Omega$.
Now, it follows from the definition of the function $F_{i,j}(\cdot,\cdot)$ that
$$F_{i,j}(\hat x,\d_1)\neq F_{i,j}(\hat x,\d_2),\;\; {\rm when}~\d_1\;\mbox{and}\;\d_2 \;\mbox{are different numbers in}\;\;(0,\d).$$
Since $(\hat x, \d_1)$ and $(\hat x, \d_2)$ belong to $D$, $F_{i,j}$ is not identically zero over $D$ for each pair $(i,j)$.

\medskip

\noindent{\it Step 3. To prove (\ref{W4.1})}.\\

Since each   $F_{i,j}$ is real analytic and is not identically a constant over $D$, the set
$$W_{i,j}\triangleq\left\{(x,\rho)\in D\bigm|F_{i,j}(x,\rho)+\lan\chi_{\omega}\xi_i,\xi_j\ran=0\right\}$$
is a real analytic subvariety with dimension at most $n$. Thus, the $\mathbb{R}^{n+1}-$Lebesgue measure of the set
$$W\triangleq\bigcup_{i,j}W_{i,j}$$
is  zero. Write
$W^\rho=\left\{x\in \Omega \bigm|(x,\rho)\in W\right\}$.
Denote by  $m(W^\rho)$ the  $\mathbb{R}^n-$Lebesgue measure  of $W^\rho$.
According to Fubini's Theorem,
$$0=\ds\int_{\om\times(0,\infty)}\chi_W(x,\rho)dxd\rho=\int^{\infty}_0\int_{\om}\chi_{W^\rho}(x)dxd\rho
=\int^{\infty}_0m(W^\rho)d\rho.$$
Thus, there is a subset $E\subset(0,\infty)$ of zero measure such that $m(W^\rho)=0$ for each $\rho\in (0,\infty)\setminus E $.
On the other hand, since $\om\setminus\overline{\omega}\neq\varnothing$,  there is $\bar\rho>0$ such that
 $\om^\rho \neq\varnothing$ for all $\rho\in (0, \bar\rho]$.

 Now, for each  $\e>0$, we arbitrarily take an  $\varepsilon_0$ from the set $(0,\min\{\bar\rho,\e\})\setminus E$.
Then, it stands that  $m\left(W^{\varepsilon_0}\right)=0$ and $\om^{\varepsilon_0}\neq \varnothing$. Hence,
\begin{eqnarray}\label{W4.4}
m\left(\Omega^{\varepsilon_0}\setminus W^{\varepsilon_0}\right)=m\left(\Omega^{\varepsilon_0}\right).
\end{eqnarray}
 Clearly, the statement that $x\in \Omega^{\varepsilon_0}\setminus W^{\varepsilon_0}$ is equivalent to the statement that
 $$
 x\in \Omega^{\varepsilon_0}\;\;\mbox{and}\;\; F_{i,j}(x,\varepsilon_0)+ \lan\chi_{\omega}\xi_i,\xi_j\ran\neq 0\;\;\mbox{for all}\;\; i, j\in \mathbb{N}.
 $$
 This, together with (\ref{W4.4}), yields that
 \begin{eqnarray}\label{W4.5} F_{i,j}(\widetilde{x},\varepsilon_0)+\lan\chi_{\omega}\xi_i,\xi_j\ran\neq0\;\;\mbox{for all}\;\; i,j\in \mathbb{N}\;\;\mbox{and for almost every}\;\; \widetilde{x}\in \om^{\varepsilon_0}.
 \end{eqnarray}

Finally, by the definition of $\om^{\varepsilon_0}$, we see that
$$B_{\varepsilon_0}(x)\subset\om\;\;\mbox{and}\;\; B_{\varepsilon_0}(x)\bigcap\omega=\varnothing\;\;\mbox{for all}\;\; x\in \om^{\varepsilon_0}.$$
Along with (\ref{W4.5}), these indicate (\ref{W4.1}). This completes the proof.
\endpf
\bigskip

\begin{Remark}
Let $\{a_i\}_{i=1}^\infty\in l^2_+\triangleq\left\{\{b_i\}_{i=1}^\infty\in l^2\; \big|\; b_i>0\;\;\mbox{for all}\;\; i\right\}$.
Consider the following time optimal control problem $(P)$: $\ds\inf\limits\left\{t: y(t; u,y_0)=0\right\}$,
 where the infimum is taken over all $u$ from the set:
 $$
 U_{ad}=\left \{u=\sum_{i=1}^\infty u_i(t)\xi_i \biggm| \;\mbox{each}\;\;u_i(\cdot)\;\;\mbox{is measurable from}\;\; \mathbb{R}^+\;\;\mbox{to}\;\;[-a_i,a_i]\right\},
 $$
 and $y(\cdot; u,y_0)$ is the solution to Equation (\ref{state}). The set $U_{ad}$ is called a control constraint set of the rectangular type. We say Problem $(P)$ has the bang-bang property if
 any optimal control $u^*=\sum_{i=1}^\infty u^*_i(t)\xi_i$ satisfies that for each $i$, $u^*_i(t)=a_i$ for a.e. $t\in (0,t^*)$, where $t^*$ is the optimal time.

{ It  is not clear to us what conditions are needed to obtain  the bang-bang property for Problem $(P)$.}
With regard to this question, we would like to mention the following: $(i)$ It is necessary to impose certain conditions on $\{a_i\}_{i=1}^\infty\in l^2_+$ to ensure the existence of optimal controls for Problem $(P)$ (\see \cite{Lv}); $(ii)$ When $U_{ad}$ is replaced by the following control constraint sets of the ball type:
$$
\widetilde{U}_{ad}\triangleq \left\{ u(\cdot)\in L^\infty(\mathbb{R}^+; L^2(\Omega)) \biggm| u(t)\in \widetilde{B}(0,r)\right\},
$$
where $\widetilde{B}(0,r)$ is the ball in $L^2(\Omega)$, centered at the origin and of radius $r>0$, the bang-bang property for the corresponding time optimal control problem $(P)$ has been studied (see  \cite{Fattorini1}, \cite{Mizel}, \cite{Wang} and \cite{Phung}).

\end{Remark}

\end{document}